\documentclass{cmslatex}
\usepackage[paperwidth=7in, paperheight=10in, margin=.875in]{geometry}
 \usepackage[backref,colorlinks,linkcolor=red,anchorcolor=green,citecolor=blue]{hyperref}
\usepackage{amsfonts,amssymb}
\usepackage{amsmath}
\usepackage{graphicx}
\usepackage{cite}
\usepackage{enumerate}
\renewcommand{\theequation}{\arabic{section}.\arabic{equation}}

   \allowdisplaybreaks
\begin{document}
 \title{On violent instability of a plasma-vacuum interface for an incompressible plasma flow  and a nonzero displacement current in vacuum\thanks{Received date, and accepted date (The correct dates will be entered by the editor).}}


          \author{Yuri Trakhinin\thanks{Sobolev Institute of Mathematics, Koptyug av. 4, 630090 Novosibirsk, Russia and  Novosibirsk State University, Pirogova str. 1, 630090 Novosibirsk, Russia (trakhinin@mail.ru).}}

         \pagestyle{myheadings} \markboth{ON VIOLENT INSTABILITY OF A PLASMA-VACUUM INTERFACE}{YURI TRAKHININ} \maketitle

          \begin{abstract}
          In the classical statement of the plasma-vacuum interface problem in ideal magnetohydrodynamics (MHD) one neglects the displacement current in the vacuum region that gives the div-curl system of pre-Maxwell dynamics for the vacuum magnetic field. For understanding the influence of the vacuum electric field on the evolution of a plasma-vacuum interface we do not neglect the displacement current and  consider the full Maxwell equations in  vacuum. For the case of an incompressible plasma flow, by constructing an Hadamard-type ill-posedness example for the constant coefficient linearized problem  we find a necessary and sufficient condition for the violent instability of a planar plasma-vacuum interface. In particular, we prove that as soon as the unperturbed plasma and vacuum magnetic fields are collinear, any nonzero unperturbed vacuum electric field makes the planar interface violently unstable. This shows the necessity of the corresponding non-collinearity condition for well-posedness and a crucial role of the vacuum electric field in the evolution of a plasma-vacuum interface.
          \end{abstract}


\section{Introduction}

\label{s1}

We consider the equations of ideal incompressible MHD, i.e., the equations governing the motion
of a perfectly conducting inviscid incompressible fluid (in particular, plasma) in magnetic field:
\begin{align}
\bar{\rho}\,  &\frac{{\rm d}v}{{\rm d}t}-(H\cdot \nabla )H+\nabla q =0, \label{1}\\[3pt]
& \frac{{\rm d}H}{{\rm d}t}-(H\cdot\nabla )v=0, \label{2}\\
& {\rm div}\,v=0,\quad {\rm div}\,H=0,\label{3}
\end{align}
where $\bar{\rho}= {\rm const} >0$ denotes density, ${\rm d}/{\rm d}t= \partial_t +v\cdot\nabla$ material derivative, $v=(v_1,v_2,v_3)$ fluid velocity, $H=(H_1,H_2,H_3)$ magnetic field,  $q =p+\frac{1}{2}|{H} |^2$ total pressure, and $p$ pressure. Note that the second equation in \eqref{3} is the {\it divergence constraint} on the initial data $U(0,x)=U_0(x)$, where $U=(v,H)$.

As is known, in the process of derivation of the nonrelativistic MHD equations one neglects the displacement current $\frac{1}{c}\partial_tE$ (see, e.g., \cite{Goed,Sed}), where ${E}=(E_1,E_2,E_3)$ is the plasma electric field and $c$ is the speed of light in vacuum. Then, $E$  is a secondary variable that can be computed from the relation
\begin{equation}
\label{4}
E=-\frac{1}{c}\,v\times H
\end{equation}
after finding $v$ and $H$ from the MHD system, and in the nonrelativistic setting $|E|^2\ll |H|^2$.

It is quite natural that in the classical statement \cite{BFKK,Goed} of the plasma-vacuum interface problem one neglects the displacement current $\frac{1}{c}\partial_t\mathcal{E}$ also in the vacuum region and instead of the Maxwell system
\begin{align}
& \frac{1}{c}\,\partial_t\mathcal{H}+\nabla\times \mathcal{E} =0, \label{5}\\[3pt]
& \frac{1}{c}\,\partial_t\mathcal{E}-\nabla\times \mathcal{H} =0, \label{6}\\
& {\rm div}\,\mathcal{H}=0,\quad {\rm div}\,\mathcal{E}=0 \label{7}
\end{align}
for the vacuum magnetic and electric fields $\mathcal{H}=(\mathcal{H}_1,\mathcal{H}_2,\mathcal{H}_3)$ and $\mathcal{E}=(\mathcal{E}_1,\mathcal{E}_2,\mathcal{E}_3)$ considers the equations
\begin{equation}
\nabla\times \mathcal{H}=0,\quad {\rm div}\,\mathcal{H}=0
\label{8}
\end{equation}
of {\it pre-Maxwell dynamics}. In this case the electric field $\mathcal{E}$ is again a second variable that can be found from \eqref{5} and the second equation in \eqref{7}. Note that equations \eqref{7} are just the divergence constraints on the initial data for system \eqref{5}, \eqref{6} whereas in the setting of pre-Maxwell dynamics they are among the basic equations.

The classical statement of the plasma-vacuum problem for systems \eqref{1}--\eqref{3} and \eqref{8} is closed by the boundary conditions
\begin{align}
& \frac{{\rm d}F }{{\rm d} t}=0,\quad [q]=0,\quad  H\cdot N =0 \label{9}\\
& \mathcal{H}\cdot  N=0, \label{10}
\end{align}
on an interface $\Gamma (t)=\{F(t,x)=0\}$ and the initial data
\begin{equation}
\label{11}
\begin{array}{ll}
{U} (0,{x})={U}_0({x}),\quad {x}\in \Omega^{+} (0),\quad
F(0,{x})=F_0({x}),\quad {x}\in\Gamma(0) , \\
\mathcal{H}(0,x)=
\mathcal{H}_0(x),\quad {x}\in \Omega^{-}(0),
\end{array}
\end{equation}
for the plasma variable $U$, the vacuum magnetic field $\mathcal{H}$ and the function $F$, where $\Omega^+(t)$ and $\Omega ^-(t)$ are space-time domains occupied by the plasma and the vacuum respectively, $N=\nabla F$, and $[q]= (q-\frac{1}{2}|\mathcal{H}|^2)|_{\Gamma}$ denotes the jump of the total pressure across the interface. The first condition in \eqref{9} means that the interface moves with the velocity of plasma particles at the boundary and since $F$ is an unknown, problem \eqref{1}--\eqref{3}, \eqref{8}--\eqref{11} is a free-boundary problem. The last condition in \eqref{9} can be shown to be the {\it boundary constraint} on the initial data \eqref{11} (see \cite{MTT1} for an analogous proof for incompressible current-vortex sheets).

Finding stability criteria of equilibrium states for a plasma-vacuum system was very popular in the 1950--70's in the context of the plasma confinement problem. The typical work in this direction is the classical paper of Bernstein et. al. \cite{BFKK}. However, the first mathematical study of the well-poseness of the ({\it non-stationary}) plasma-vacuum interface problem was carried out relatively recently in \cite{T10} for compressible MHD. In \cite{T10} a basic energy a priori estimate in Sobolev spaces for the linearized plasma-vacuum interface problem was proved under the non-collinearity condition
\begin{equation}\label{noncollin}
|{H}\times{\mathcal{H}}| \geq \epsilon_1>0\quad \mbox{on}\ \Gamma (t)
\end{equation}
satisfied on the interface for the unperturbed flow (by the unperturbed flow we mean a basic state about which one linearizes the problem). In \cite{T10} the alternative well-posedness condition
\begin{equation}\label{RT}
\left[\frac{\partial q}{\partial N}\right] \leq -\epsilon_2<0\quad \mbox{on}\ \Gamma (t),
\end{equation}
which is nothing else than the (generalized) Rayleigh-Taylor sign condition, was also proposed, where $[\partial q/\partial N] = \big(\partial q/\partial N- \frac{1}{2} \partial|\mathcal{H}|^2/\partial N\big)|_{\Gamma}$.

Under the non-collinearity condition \eqref{noncollin} satisfied at the first moment, the well-posedness of the nonlinear problem (for compressible MHD) was proved by Secchi and Trakhinin \cite{ST2} basing on their preparatory well-posedness result \cite{ST1} for the linearized problem. The proof of the local-time well-posedness of the plasma-vacuum interface problem under the Rayleigh-Taylor sign condition \eqref{RT} is still an open problem not only for compressible MHD but also for an incompressible plasma flow (for compressible MHD, some particular results in this direction can be found in \cite{T10,Tcpaa}).
In comparison with the study of the compressible plasma-vacuum problem \cite{ST1,ST2,T10,Tcpaa}, much more attention was paid by researches to the case of incompressible MHD for which the plasma-vacuum system is modelled by problem \eqref{5}--\eqref{11}. In this connection, we can refer the reader to results for this problem obtained in \cite{Gu1,Gu2,GuWang,Hao,HaoLuo,HaoLuo2,MTTpv,SWZ2}. In our opinion, the most complete result was recently obtained by Sun, Wang and Zhang \cite{SWZ2} who proved the well-posedness of the incompressible plasma-vacuum problem, provided that the initial data satisfy the non-collinearity condition \eqref{noncollin}.

The more general ``stability'' assumption proposed in \cite{Tcpaa} requires that the Rayleigh-Taylor sign condition \eqref{RT} is satisfied only at all those points of the initial interface where the non-collinearity condition \eqref{noncollin} fails. The proof of the local well-posedness of the plasma-vacuum interface problem under this assumption remains an open problem for both compressible and incompressible MHD. By constructing an Hadamard-type ill-posedness example for the frozen coefficients linearized problem (for both compressible and incompressible MHD) it was  proved in \cite{Tcpaa} that the simultaneous  failure of the non-collinearity condition and the Rayleigh-Taylor sign condition leads to Rayleigh-Taylor instability, more precisely, ill-posedness takes place if and only if both of these conditions are violated for frozen coefficients. In the light of this, it is natural to suppose that the hypothetical well-posedness result under the more general ``stability'' assumption will be one day obtained.

On the other hand, the main goal of the present paper is to show that from the physical point of view the non-collinearity condition \eqref{noncollin} is not only sufficient but also {\it necessary for well-posedness}. For this purpose, we consider the technically simpler case of incompressible MHD for which, unlike the compressible case studied by Mandrik and Trakhinin \cite{MT} being motivated by results in \cite{T12} for the relativistic setting, it is possible to perform a complete spectral analysis of the linearized problem with constant coefficients, provided that the displacement current $\frac{1}{c}\partial_t\mathcal{E}$ was not neglected in the vacuum region and in the nonlinear problem the vacuum magnetic and electric fields $\mathcal{H}$ and $\mathcal{E}$  satisfy the Maxwell system \eqref{5}--\eqref{7}.

We note that the same compressible plasma-vacuum problem with a nonzero displacement current in vacuum as in \cite{MT} (see also \cite{M1}) was independently studied by Catania, D'Abbicco and Secchi \cite{CDS,CDS2}. In this paper we consider the counterpart of the plasma-vacuum problem from \cite{CDS,CDS2,M1,MT}  for incompressible MHD. Actually, this problem for incompressible MHD (with a nonzero displacement current in vacuum) was already being considered by Secchi \cite{Sec_sw} for the case of two space dimensions in the context of the study of weakly nonlinear surface waves on a plasma-vacuum interface. However, for the two-dimensional (2D) case the normal component (with respect to the interface) of the vacuum electric field playing the crucial role
in the appearance of violent instability/ill-posedness in the three-dimensional (3D) case (found in \cite{MT} for compressible MHD) is zero by definition. In this sense the 2D case is exceptional whereas in this paper we study the general 3D case and show, in particular, that as soon as the unperturbed plasma and vacuum magnetic fields are collinear, any nonzero unperturbed vacuum electric field makes the planar interface violently unstable. This shows the necessity of the non-collinearity condition for well-posedness and a crucial role of the vacuum electric field in the evolution of a plasma-vacuum interface. Moreover, we find a necessary and sufficient condition (on the normal component of the vacuum electric field) for the violent instability of a planar plasma-vacuum interface.

\section{Statement of the plasma-vacuum problem for a nonzero displacement current in vacuum}

Following \cite{CDS,MT,Sec_sw}, we do not neglect the displacement current in the vacuum region and in $\Omega^-(t)$ consider the Maxwell equations \eqref{5}--\eqref{7}. Since in this paper we study only the linearized problem, we do not care about a proper geometry of reference domains $\Omega^{\pm}(t)$ necessary for the subsequent analysis of the original nonlinear problem. In particular, unlike, for example, \cite{SWZ2}, we do not introduce periodical boundary conditions in the tangential $x_{2,3}$--directions and rigid wall boundary conditions whose introduction for the plasma region is motivated by the fact that the total pressure $q$ is an ``elliptic'' unknown for which it is reasonable to have a problem in a bounded domain. As in \cite{CDS,MT,Sec_sw}, for technical simplicity we assume that the interface $\Gamma (t)$ has the form of a graph and the domains $\Omega^{\pm}(t)$ are unbounded:
\[
\Gamma (t)=\{F(t,x)=x_1-\varphi (t,x')\},\quad \Omega^\pm(t)=\{\pm (x_1- \varphi(t,x'))>0,\ x'\in \mathbb{R}^2\},\quad x'=(x_2,x_3).
\]
It is reasonable to reduce the MHD system and the Maxwell equations to a dimensionless form.
In terms of the scaled values
\begin{equation}
\label{12}
\tilde{x}=\frac{x}{\ell},\quad \tilde{t}=\frac{\bar{v}t}{\ell},\quad \tilde{v}=\frac{v}{\bar{v}},\quad \tilde{q}=\frac{q}{\bar{\rho}\bar{v}^2},\quad
 \widetilde{H}=\frac{H}{\bar{v}\sqrt{\bar{\rho}}},\quad \widetilde{\mathcal{H}}=\frac{\mathcal{H}}{\bar{v}\sqrt{\bar{\rho}}},\quad
\widetilde{\mathcal{E}}=\frac{\mathcal{E}}{\bar{v}\sqrt{\bar{\rho}}}
\end{equation}
and after dropping tildes systems \eqref{1}, \eqref{2} and \eqref{5}, \eqref{6} read:
\begin{equation}
\left\{
\begin{array}{ll}
{\rm div}\, v=0, & \\[6pt]
\displaystyle \frac{{\rm d}v}{{\rm d}t}-(H\cdot \nabla )H+\nabla q =0, & \\[12pt]
\displaystyle \frac{{\rm d}H}{{\rm d}t}-(H\cdot\nabla )v=0  &\quad \mbox{in}\ \Omega^+(t),
 \end{array}\right.
\label{13}
\end{equation}
\begin{equation}
\left\{
\begin{array}{l}
\varepsilon\,\partial_t\mathcal{H}+\nabla\times \mathcal{E} =0, \\[3pt]
\varepsilon\,\partial_t\mathcal{E}-\nabla\times \mathcal{H} =0\qquad  \mbox{in}\ \Omega^-(t),
 \end{array}\right.
\label{14}
\end{equation}
where $\ell = {\rm const}>0 $ is a characteristic length, $\bar{v} = {\rm const}>0 $ is a characteristic (average) speed of the plasma flow, and $\varepsilon =\bar{v}/c \ll 1$ is the natural small parameter of the problem.

The boundary conditions (in terms of the scaled values \eqref{12}) are the same as for the case of compressible MHD in \cite{CDS,MT}:
\begin{align}
& \partial_t \varphi= v_N,\label{15'}\\ & q= \textstyle{\frac{1}{2}}\left(|\mathcal{H}|^2-|\mathcal{E}|^2\right),   \label{15}\\
& \mathcal{E}_{\tau_2}=\varepsilon\mathcal{H}_3\partial_t\varphi,\quad \mathcal{E}_{\tau_3}=-\varepsilon\mathcal{H}_2\partial_t\varphi \qquad \mbox{on}\ \Gamma (t),\label{16}
\end{align}
where $v_N=v\cdot N$, $N=\nabla F= (1,-\partial_2\varphi ,-\partial_3\varphi )$ and $\mathcal{E}_{\tau_i}=\mathcal{E}_1\partial_i\varphi+\mathcal{E}_i$ ($i=2,3$).
The statement of the plasma-vacuum interface problem is closed by  the initial data
\begin{equation}
\label{17}
\begin{array}{ll}
{U} (0,{x})={U}_0({x}),\quad {x}\in \Omega^{+} (0),\quad
\varphi (0,x')=\varphi_0(x'),\quad x'\in\mathbb{R}^2 , \\
V(0,x)=
V_0(x),\quad {x}\in \Omega^{-}(0)
\end{array}
\end{equation}
for $U=(v,H)$, $V=(\mathcal{H},\mathcal{E})$ and the function $\varphi$. At last, as for the case of compressible plasma flow in \cite{CDS,MT}, one can show that
\begin{equation}
{\rm div}\,{H}=0\quad \mbox{in}\ \Omega^+(t),\quad {\rm div}\,\mathcal{H}=0,\quad {\rm div}\,\mathcal{E}=0\quad \mbox{in}\ \Omega^-(t)
\label{18}
\end{equation}
and
\begin{equation}
{H}_N=0,\quad \mathcal{H}_N=0\quad \mbox{on}\ \Gamma (t)
\label{19}
\end{equation}
are the divergence and boundary constraints on the initial data \eqref{17}, i.e., they hold for $t>0$ if they were satisfied at $t=0$. Here
${H}_N=H\cdot N$ and $\mathcal{H}_N=\mathcal{H}\cdot N$. In fact, problem \eqref{13}--\eqref{17} was also written down in \cite{Sec_sw} where the analysis of the amplitude equation for surface waves was then performed for the 2D case.

Note the boundary condition \eqref{15} expresses the fact that there is no jump of the total pressure across the interface, where ${\frac{1}{2}}(|\mathcal{H}|^2-|\mathcal{E}|^2)$ can be interpreted as the magnetic pressure in vacuum (see Remark \ref{r1} below). The boundary conditions \eqref{16} follow from the jump conditions \cite{BFKK,Sed}
\begin{equation}
N\times [E]=\varepsilon\partial_t\varphi\,[H] \quad \mbox{on}\ \Gamma (t)
\label{21}
\end{equation}
for the conservation laws
\[
\varepsilon\partial_t H^{\pm} +\nabla \times E^{\pm} =0\quad \mbox{in}\ \Omega^{\pm}(t)
\]
if we take into account \eqref{15'} and the first boundary constraint in \eqref{19}, where $H^+=H$, $H^-=\mathcal{H}$, $E^+= E=-\varepsilon (v\times H)$ (cf. \eqref{4}), $E^-=\mathcal{E}$, and the above conservation laws in $\Omega^+(t)$ are just the conservative form of the last three scalar equations of system \eqref{13} obtained with the help of the first divergence constraint in \eqref{18}. Moreover, in view of the first constraint in \eqref{19}, the first condition in \eqref{21} is nothing else than the second boundary constraint in \eqref{19}.

\begin{remark}
{\rm In relativistic MHD the total pressure in plasma reads \cite{Goed,T12}:
\[
q=p +p_m = p +\frac{|H|^2}{2}\left(1-\frac{|v|^2}{c^2}\right) +\frac{1}{2}\left(\frac{v}{c}\cdot H\right)^2 = p +\frac{1}{2}\left(H^2-E^2\right),
\]
where $E$ is given by \eqref{4}. In vacuum  the value  corresponding to $p_m$ is $q_v=\frac{1}{2}(|\mathcal{H}|^2-|\mathcal{E}|^2)$. In the nonrelativistic limit $p_m=\frac{1}{2}|H|^2$ whereas $q_v$ stays the same as in the relativistic setting. Note also that $|\mathcal{H}|^2-|\mathcal{E}|^2$ is one of the two fundamental Lorentz invariants of the electromagnetic field (the second one is $\mathcal{H}\cdot\mathcal{E}$; see, e.g., \cite{Sed}).
}
\label{r1}
\end{remark}

The rest of the paper is organized as follows. In the next section we write down the linearization of problem \eqref{13}--\eqref{17} and formulate the main result about its ill-posedness for the case of constant coefficients (see Theorem \ref{t1}). In Sect. \ref{s4}, by spectral analysis we prove Theorem \ref{t1}. At last, in Sect. \ref{s5} we briefly discuss open problems connected with the proof of well-posedness of the linearized problem and the nonlinear problem for the vacuum electric field satisfying the hypothetical ``stability'' condition \eqref{well_nonlin}.

\section{Linearized problem and main result}

\label{s3}

For our present goals it is enough to straighten the interface $\Gamma$ by using the simplest change of independent variables
\begin{equation}
\tilde{t}=t,\quad \tilde{x}_1= x_1-\varphi (t,x'),\quad \tilde{x}'=x'
\label{chv}
\end{equation}
(more involved changes of independent variables can be found, for example, in \cite{ST2,SWZ1,SWZ2,Tcpam}). By this change we reduce our free boundary problem \eqref{13}--\eqref{17} to the following initial-boundary value problem in the fixed domains $[0,T]\times \mathbb{R}^3_{\pm}$ with $\mathbb{R}^3_{\pm}=\{\pm x_1>0, x'\in\mathbb{R}^2\}$ (here and below we drop tildes):
\begin{equation}
\left\{ \begin{array}{l}
{\rm div}\,w=0,\\
\partial_tU +A_n(U,\varphi)\partial_1U+A_2(U)\partial_2U+A_3(U)\partial_3U +\nabla_Nq=0\quad \mbox{in}\ [0,T]\times \mathbb{R}^3_+,
\end{array}\right.
\label{22}
\end{equation}
\begin{equation}
\varepsilon \partial_tV +B_n(\varphi)\partial_1V+B_2\partial_2V+B_3\partial_3V =0\quad \mbox{in}\ [0,T]\times \mathbb{R}^3_-,
 \label{23}
\end{equation}
\begin{equation}
\left\{ \begin{array}{l}
\partial_t \varphi= v_N,\quad q= \textstyle{\frac{1}{2}}\left(|\mathcal{H}|^2-|\mathcal{E}|^2\right),   \\
\mathcal{E}_{\tau_2}=\varepsilon\mathcal{H}_3\partial_t\varphi,\quad \mathcal{E}_{\tau_3}=-\varepsilon\mathcal{H}_2\partial_t\varphi \qquad\mbox{on}\ [0,T]\times \{x_1=0\}\times\mathbb{R}^2,
\end{array}\right.
\label{24}
\end{equation}
\begin{equation}
\label{25}
\begin{array}{ll}
{U} (0,{x})={U}_0({x}),\quad {x}\in \mathbb{R}^3_{+},\quad
\varphi (0,x')=\varphi_0(x'),\quad x'\in\mathbb{R}^2 , \\
V(0,x)=
V_0(x),\quad {x}\in \mathbb{R}^3_-,
\end{array}
\end{equation}
where $w=(v_N,v_2,v_3)$, $\nabla_Nq= (\partial_1q)N + (0,\partial_2q,\partial_3q)$,
\[
A_k(U)=\left( \begin{array}{cc} v_k & -H_k\\
-H_k & v_k\\
\end{array} \right)\otimes I_3,\quad k=1,2,3,\quad
B_1=\left(\begin{array}{cccccc}
0 & 0 & 0& 0 & 0 & 0 \\
0 & 0 & 0& 0 & 0 & -1 \\
0 & 0 & 0& 0 & 1 & 0 \\
0 & 0 & 0& 0 & 0 & 0 \\
0 & 0 & 1& 0 & 0 & 0 \\
0 & -1 & 0& 0 & 0 & 0
\end{array} \right),
\]
\[
B_2=\left(\begin{array}{cccccc}
0 & 0 & 0& 0 & 0 & 1 \\
0 & 0 & 0& 0 & 0 & 0 \\
0 & 0 & 0& -1 & 0 & 0 \\
0 & 0 & -1& 0 & 0 & 0 \\
0 & 0 & 0& 0 & 0 & 0 \\
1 & 0 & 0& 0 & 0 & 0
\end{array} \right),\quad
B_3=\left(\begin{array}{cccccc}
0 & 0 & 0& 0 & -1 & 0 \\
0 & 0 & 0& 1 & 0 & 0 \\
0 & 0 & 0& 0 & 0 & 0 \\
0 & 1 & 0& 0 & 0 & 0 \\
-1 & 0 & 0& 0 & 0 & 0 \\
0 & 0 & 0& 0 & 0 & 0
\end{array} \right),
\]
\[
A_n(U,\varphi )= A_1(U)- I_6 \partial_t\varphi -A_2(U) \partial_2\varphi-A_3(U)\partial_3\varphi,
\]
\[
B_n(\varphi )= B_1(U)- \varepsilon I_6 \partial_t\varphi -B_2\partial_2\varphi-B_3\partial_3\varphi,
\]
and $I_j$ is the unit matrix of order $j$. Moreover, for the initial data \eqref{25} we have the boun\-dary constraints \eqref{19} at $x_1=0$ and the divergence constraints written for the ``curved'' fields  $(H_N,H_2,H_3)$, $(\mathcal{H}_N,\mathcal{H}_2,\mathcal{H}_3)$ and $(\mathcal{E}_N,\mathcal{E}_2,\mathcal{E}_3)$.

The first step in the study of the well-posedness of a nonlinear initial boundary value problem is the spectral analysis of the corresponding constant coefficient linearized problem for finding a possibly non-empty domain of parameters in which this linear problem is ill-posed. The constant coefficient problem results from the linearization of the nonlinear problem about its constant solution. If the original nonlinear problem is a free boundary/interface problem, then this constant solution is usually associated with a planar interface. If the interface is described by the equation $x_1=\varphi (t,x')$ and if the interior equations in the domains $x_1\gtrless \varphi (t,x')$ are Galilean invariant, then, without loss of generality, one can consider the constant solution associated with the planar interface $x_1=0$.

Regarding our plasma-vacuum interface problem, the Maxwell equations {\it are not} Galilean invariant (they are Lorentz invariant). That is, we have to consider the planar interface $x_1=\sigma t$ moving with a constant speed $\sigma$. The corresponding constant solution $(U,V)=(\widehat{U},\widehat{V})$ (with $\varphi =\sigma t$) of the reduced nonlinear problem \eqref{22}--\eqref{24} (if we also take into account the boundary constraints) reads:
\[
\widehat{U}= (\hat{v},\widehat{H}),\quad \widehat{V}= (\widehat{\mathcal{H}},\widehat{\mathcal{E}}),
\]
with
\begin{equation}
\begin{split}
& \hat{v}=(\sigma ,\hat{v}'),\quad \widehat{H}=(0,\widehat{H}'),\quad \widehat{\mathcal{H}}=(0,\widehat{\mathcal{H}}'),\quad
\widehat{\mathcal{E}}=(\widehat{\mathcal{E}}_1,\varepsilon\sigma\widehat{\mathcal{H}}_3,-\varepsilon\sigma
\widehat{\mathcal{H}}_2),
\\
& \hat{v}'=(\hat{v}_2,\hat{v}_3),\quad \widehat{H}'=(\widehat{H}_2,\widehat{H}_3),\quad \widehat{\mathcal{H}}'=(\widehat{\mathcal{H}}_2,\widehat{\mathcal{H}}_3),
\end{split}
\label{constsol}
\end{equation}
where $\hat{v}_k$, $\widehat{H}_k$, $\widehat{\mathcal{H}}_k$ ($k=2,3$) and $\widehat{\mathcal{E}}_1$ are some constants. Moreover, as follows from the physical condition $p_{|x_1=0}\geq 0$ and the second boundary condition in \eqref{24}, we have the restriction
\[
|\widehat{\mathcal{H}}|^2(1-\varepsilon^2\sigma^2)\geq\widehat{\mathcal{E}}_1^{\,2}+|\widehat{H}|^2.
\]
Hence, it is necessary that
\begin{equation}
|\widehat{\mathcal{H}}|^2\geq\widehat{\mathcal{E}}_1^{\,2}.
\label{physrestr}
\end{equation}
We will below show (see Remark \ref{r4}) that as soon as the physical restriction \eqref{physrestr} is violated the corresponding linearized constant coefficient problem is ill-posed.

\begin{remark}
{\rm
It is interesting to note that in \cite{MT} condition \eqref{physrestr}  was not taken into account for the compressible plasma-vacuum interface problem as a physical restriction, but a part of the ill-posedness domain of the linearized problem found there for some particular cases of the constant solution is described exactly by the inequality opposite to \eqref{physrestr}.
}
\label{r0}
\end{remark}

Linearizing problem \eqref{22}--\eqref{24} about its exact solution $(\widehat{U},\widehat{V},\sigma t)$, we obtain the following constant coefficient problem for the perturbations $U$, $V$ and $\varphi$ which are denoted by the same letters as the unknowns of the nonlinear problem:
\begin{equation}
\left\{ \begin{array}{l}
{\rm div}\,v=0,\\
\partial_tU +\widehat{A}_2\partial_2U+\widehat{A}_3\partial_3U +\nabla q=0\quad \mbox{in}\ [0,T]\times \mathbb{R}^3_+,
\end{array}\right.
\label{26}
\end{equation}
\begin{equation}
\varepsilon (\partial_t-\sigma\partial_1)V +\sum_{j=1}^{3}B_j\partial_jV =0\quad \mbox{in}\ [0,T]\times \mathbb{R}^3_-,
 \label{27}
\end{equation}
\begin{equation}
\left\{ \begin{array}{ll}
\partial_t \varphi= v_1-\hat{v}_2\partial_2\varphi -\hat{v}_3\partial_3\varphi, & \\
q= \widehat{\mathcal{H}}_2(\mathcal{H}_2+\varepsilon \sigma \mathcal{E}_3) + \widehat{\mathcal{H}}_3(\mathcal{H}_3-\varepsilon \sigma \mathcal{E}_2) -\widehat{\mathcal{E}}_1\mathcal{E}_1,  & \\
\widehat{\mathcal{E}}_1\partial_2\varphi +\mathcal{E}_2-\varepsilon \sigma \mathcal{H}_3=\varepsilon\widehat{\mathcal{H}}_3\partial_t\varphi, & \\ \widehat{\mathcal{E}}_1\partial_3\varphi +\mathcal{E}_3+\varepsilon \sigma \mathcal{H}_2=-\varepsilon\widehat{\mathcal{H}}_2\partial_t\varphi & \mbox{on}\ [0,T]\times\{x_1=0\}\times \mathbb{R}^2,
\end{array}\right.
\label{28}
\end{equation}
with corresponding initial data in the same form as in \eqref{25}, where $\widehat{A}_k=A_k (\widehat{U})$ ($k=2,3$). Moreover, one can show that
\begin{equation}
{\rm div}\,{H}=0\quad \mbox{in}\ [0,T]\times \mathbb{R}^3_+,\quad {\rm div}\,\mathcal{H}=0,\quad {\rm div}\,\mathcal{E}=0\quad \mbox{in}\ [0,T]\times \mathbb{R}^3_-
\label{18l}
\end{equation}
and
\begin{equation}
{H}_1=\widehat{H}_2\partial_2\varphi+\widehat{H}_3\partial_3\varphi,\quad \mathcal{H}_1=\widehat{\mathcal{H}}_2\partial_2\varphi+\widehat{\mathcal{H}}_3\partial_3\varphi\quad \mbox{on}\ [0,T]\times\{x_1=0\}\times \mathbb{R}^2
\label{19l}
\end{equation}
hold for $t>0$ if they were satisfied for the initial data of problem \eqref{26}--\eqref{28}.

Following \cite{CDS}, we can introduce the new unknowns
\begin{equation}
\begin{split}
\breve{\mathcal{H}} & =(\mathcal{H}_1,\mathcal{H}_2+\varepsilon\sigma \mathcal{E}_3,\mathcal{H}_3-\varepsilon\sigma \mathcal{E}_2),\\ \breve{\mathcal{E}} & = (\mathcal{E}_1,\mathcal{E}_2-\varepsilon\sigma \mathcal{H}_3,\mathcal{E}_3+\varepsilon\sigma\mathcal{H}_2).
\end{split}
\label{30}
\end{equation}
This is nothing else than the usage of the nonrelativistic version of the Joules-Bernoulli equations (see, e.g., \cite{Sed}). In fact, in \cite{CDS} a more involved (``curved'') variant of \eqref{30} was utilized   for showing that the counterpart of our plasma-vacuum interface problem for compressible MHD has a correct number of boundary conditions. In fact, the arguments in \cite{CDS} take also place for our case of incompressible MHD, and the number of boundary conditions in \eqref{24} and \eqref{28} is correct regardless of the sign of the interface speed (${\rm sign}\,\partial_t\varphi$ for the nonlinear problem and ${\rm sign}\,\sigma$ for the linear one).

After making the change of unknowns \eqref{30} system \eqref{26} stays, of course, unchanged whereas the boundary conditions \eqref{28} after dropping breves coincide with their form for $\sigma =0$. The Maxwell equations \eqref{27} can be written as
\[
\varepsilon B_0\partial_t\breve{V} +\sum_{j=1}^{3}B_j\partial_j\breve{V} =0,
\]
where $\breve{V}=(\breve{\mathcal{H}},\breve{\mathcal{E}})$ and the matrix
\[
B_0=\frac{1}{1-\varepsilon^2\sigma^2}
\begin{pmatrix}
1-\varepsilon^2\sigma^2 & 0 & 0 & 0 & 0 & 0 \\
0 & 1 & 0 & 0 & 0 & -\varepsilon\sigma \\
0 & 0 & 1 & 0 & \varepsilon\sigma & 0 \\
0 & 0 & 0 & 1-\varepsilon^2\sigma^2 & 0 & 0\\
0 & 0& \varepsilon\sigma & 0 & 1 & 0\\
0 & -\varepsilon\sigma & 0 & 0 & 0& 1
\end{pmatrix}
\]
is found from the relation $V=B_0\breve{V}$. The matrix $B_0$ is symmetric and since $B_0|_{\varepsilon =0}=I_6$, in the nonrelativistic limit we have $B_0>0$, i.e., the system is symmetric hyperbolic. By default we, of course,  consider the nonrelativistic limit $\varepsilon \rightarrow 0$ and, as was already mentioned above, $\varepsilon $ is our natural small parameter. In particular, this means that $\varepsilon |\sigma| \ll 1$. In the nonrelativistic limit spectral properties of the above constant coefficient hyperbolic system coincide for $\sigma =0$ and $\sigma \neq 0$ (for nonrelativistic speeds $\sigma$).

Without loss of generality we may thus assume that $\sigma =0$. In fact, after finding a necessary and sufficient ill-posedness condition for problem \eqref{26}--\eqref{28} with $\sigma =0$ one can easily check that we obtain the same result for $\sigma \neq 0$. For further convenience we write down problem \eqref{26}--\eqref{28} for $\sigma =0$:
\begin{equation}
\left\{
\begin{array}{ll}
{\rm div}\, v=0, & \\[6pt]
\displaystyle \frac{{\rm d}'v}{{\rm d}t}-(\widehat{H}'\cdot \nabla' )H+\nabla q =0, & \\[12pt]
\displaystyle \frac{{\rm d}'H}{{\rm d}t}-(\widehat{H}'\cdot\nabla' )v=0  &\quad \mbox{in}\ [0,T]\times \mathbb{R}^3_+,
 \end{array}\right.
\label{33}
\end{equation}
\begin{equation}
\left\{
\begin{array}{l}
\varepsilon\,\partial_t\mathcal{H}+\nabla\times \mathcal{E} =0, \\[3pt]
\varepsilon\,\partial_t\mathcal{E}-\nabla\times \mathcal{H} =0\qquad  \mbox{in}\ [0,T]\times \mathbb{R}^3_-,
 \end{array}\right.
\label{34}
\end{equation}
\begin{equation}
\left\{ \begin{array}{ll}
 \displaystyle \frac{{\rm d}'\varphi }{{\rm d}t}= v_1, \quad q= \widehat{\mathcal{H}}_2\mathcal{H}_2+ \widehat{\mathcal{H}}_3\mathcal{H}_3 -\widehat{\mathcal{E}}_1\mathcal{E}_1,   \\[6pt]
\mathcal{E}_2=\varepsilon\widehat{\mathcal{H}}_3\partial_t\varphi - \widehat{\mathcal{E}}_1\partial_2\varphi , \quad \mathcal{E}_3=-\varepsilon\widehat{\mathcal{H}}_2\partial_t\varphi -\widehat{\mathcal{E}}_1\partial_3\varphi  \qquad \mbox{on}\ [0,T]\times\{x_1=0\}\times \mathbb{R}^2,
\end{array}\right.
\label{35}
\end{equation}
with some initial data at $t=0$, where ${\rm d}'/{\rm d}t =\partial_t +(\hat{v}'\cdot\nabla' )$ and $\nabla ' =(\partial_2,\partial_3)$. For the initial data we again have the divergence and boundary constraints \eqref{18l} and \eqref{19l}.

We are now in a position to formulate our main result which is a necessary and sufficient ill-posedness condition for problem \eqref{33}--\eqref{35}. In the below theorem violent instability means the ill-posedness of this problem whereas neutral stability means that for problem \eqref{33}--\eqref{35} the Lopatinski condition holds but the uniform Lopatinski condition \cite{Kreiss} is violated.

\begin{theorem}
For every given planar plasma-vacuum interface (described by the constants in \eqref{constsol}), there is a small value $\varepsilon_*>0$ such that  for all $\varepsilon <\varepsilon_*$ the interface is violently unstable  if and only if
\begin{equation}\label{illcond}
\widehat{\mathcal{E}}_1^{\,2}> \frac{|\widehat{H}|^2 +|\widehat{\mathcal{H}}|^2-\sqrt{\big(|\widehat{H}|^2 +|\widehat{\mathcal{H}}|^2\big)^2-4|\widehat{H}\times\widehat{\mathcal{H}}|^2}}{2}
\end{equation}
or
\begin{equation}\label{illcond=}
\widehat{\mathcal{E}}_1^{\,2}= \frac{|\widehat{H}|^2 +|\widehat{\mathcal{H}}|^2-\sqrt{\big(|\widehat{H}|^2 +|\widehat{\mathcal{H}}|^2\big)^2-4|\widehat{H}\times\widehat{\mathcal{H}}|^2}}{2}
\end{equation}
and
\begin{equation}\label{illcond=1}
\widehat{\mathcal{E}}_1 \big( (\widehat{H}\cdot\widehat{\mathcal{H}})(\hat{v}\times\widehat{H})+(|\widehat{\mathcal{H}}|^2-\widehat{\mathcal{E}}_1^{\,2})
(\hat{v}\times\widehat{\mathcal{H}})\big)\cdot e_1 >0,
\end{equation}
where $e_1=(1,0,0)$. Otherwise, for all $\varepsilon <\varepsilon_*$ this interface is neutrally stable.
\label{t1}
\end{theorem}

\begin{remark}
{\rm
In the above theorem the assertion about the existence of a small value $\varepsilon_*$ just means that our instability/stability criteria are valid in the {\it nonrelativistic limit} $\varepsilon \rightarrow 0$. However, it is technically impossible (and not physically interesting) to estimate exactly how small should be $\varepsilon_*$. We also note that our below proof of Theorem \ref{t1} is based on a perturbation argument for the small parameter $\varepsilon$. The same standard argument was, for example, used in \cite{BThand,MZ} for showing the uniform stability of fast MHD shock waves under the gas dynamical uniform stability condition in the limit $H \rightarrow 0$ of a weak magnetic field.
}
\label{r3*}
\end{remark}

As follows from the ill-posedness condition \eqref{illcond}, as soon as the unperturbed plasma and vacuum magnetic fields are collinear ($\widehat{H}\times\widehat{\mathcal{H}}=0$), any non-zero unperturbed vacuum electric field makes the planar interface violently unstable. This shows a crucial role of the vacuum electric field in the evolution of a plasma-vacuum interface.

\section{Spectral analysis of the linear problem}

\label{s4}

Instead of the classical arguments \cite{Kreiss} toward the check of the Lopatinski condition and the uniform Lopatinski condition we perform similar spectral analysis connected with the construction of an Hadamard-type ill-posedness example. Let us seek the following sequence of exponential solutions of problem \eqref{33}--\eqref{35}:
\begin{equation}
\begin{pmatrix} U_n \\ q_n \end{pmatrix} =\begin{pmatrix} \bar{U}_n \\ \bar{q}_n \end{pmatrix}\exp \left\{ n\left(s t +\lambda^+ x_1+ i(\omega ', x') \right)\right\}\quad\mbox{for}\ x_1>0,
\label{exp_sol_pl}
\end{equation}
\begin{equation}
V_n  = \bar{V}_n \exp \left\{ n\left(s t +\lambda^- x_1+ i(\omega ', x') \right)\right\}\quad\mbox{for}\ x_1<0,
\label{exp_sol_vac}
\end{equation}
\begin{equation}
\varphi_n=\bar{\varphi}_n\exp \left\{ n\left(s t + i(\omega ', x') \right)\right\},
\label{exp_sol_fr}
\end{equation}
with
\begin{equation}
\Re\,s>0,\quad \Re\,\lambda^+\leq 0,\quad \Re\,\lambda^-\geq 0,\label{re_re}
\end{equation}
where $s$, $\lambda^+$ and $\lambda^-$ are complex constants, $\omega '=(\omega_2,\omega_3)$ and $\omega_{2,3}$ are real constants,
\[
\bar{U}_n=(\bar{v}_n,\bar{H}_n),\quad \bar{V}_n=(\bar{\mathcal{H}}_n,\bar{\mathcal{E}}_n),\quad \bar{v}_n=(\bar{v}_{1n},\bar{v}_{2n},\bar{v}_{3n}),\quad \bar{H}_n=(\bar{H}_{1n},\bar{H}_{2n},\bar{H}_{3n}),
\]
\[
\bar{\mathcal{H}}_n=(\bar{\mathcal{H}}_{1n},\bar{\mathcal{H}}_{2n},\bar{\mathcal{H}}_{3n}),\quad
\bar{\mathcal{E}}_n=(\bar{\mathcal{E}}_{1n},\bar{\mathcal{E}}_{2n},\bar{\mathcal{E}}_{3n}),
\]
and $\bar{v}_{jn}$, $\bar{H}_{jn}$, $\bar{q}_n$, $\bar{\mathcal{H}}_{jn}$, $\bar{\mathcal{E}}_{jn}$, $\bar{\varphi}_n$ are complex constants ($j=1,2,3$, $n\in\mathbb{N}$). Clearly, requirements \eqref{re_re} imply the boundedness of solutions at $t=0$ and their infinite growth as $n\rightarrow \infty$ for any (even very small) $t>0$. This is nothing else than ill-posedness.

The substitution of \eqref{exp_sol_pl} into system \eqref{33} gives a dispersion relation for $s$, $\lambda^+$ and $\omega'$ which has the two roots $\lambda_{1,2}^+=\pm |\omega '|$. We could equivalently get them from the dispersion relation $(\lambda^+)^2-|\omega'|^2=0$ for the Laplace equation $\triangle q=0$ following from system \eqref{33}. That is, we choose $\lambda^+=-|\omega '|$ satisfying \eqref{re_re}. One can easily show that for problem \eqref{33}--\eqref{35} we are not able to construct a 1D Hadamard-type ill-posedness example, i.e., the sequence of exponential solutions \eqref{exp_sol_pl}--\eqref{exp_sol_fr} with $\omega'=0$ obeying requirements \eqref{re_re}. Hence, we will assume that $\omega'\neq 0$.

Substituting \eqref{exp_sol_vac} into system \eqref{34} gives the dispersion relation
\[
\det (\varepsilon sI_6 +B_1\lambda^-+i\omega_2B_2+i\omega_3B_3)=\varepsilon^2s^2(\varepsilon^2s^2 -(\lambda^-)^2)(\varepsilon^2s^2 +|\omega'|^2-(\lambda^-)^2)=0
\]
which has two roots
\[
\lambda_1^-=\varepsilon s,\quad \mbox{and}\quad
\lambda_2^-=\sqrt{|\omega '|^2+\varepsilon^2s^2}
\]
satisfying \eqref{re_re}. For the root $\lambda_1^-$ the algebraic system following from \eqref{34} formally has a non-zero solution  $\bar{V}_n$. But, if we take into account the last two necessary divergence constraints  in \eqref{18l}, then $\bar{V}_n=0$. Indeed, for $\lambda^-=\lambda^-_1$ the mentioned algebraic system implies $\bar{\mathcal{H}}_{1n}=\bar{\mathcal{E}}_{1n}$. Then, in view of the last two divergence constraints  in \eqref{18l}, at least one of the rest components of the vector $\bar{V}_{n}$ is not zero if and only if $\omega'=0$. We thus have the following roots $\lambda^{\pm}$:
\begin{equation}\label{48}
 \lambda^+= -|\omega '|,\quad \lambda^-=\sqrt{|\omega '|^2+\varepsilon^2s^2}.
\end{equation}

By substituting \eqref{exp_sol_pl}--\eqref{exp_sol_fr} into the boundary conditions \eqref{35}, one gets
\begin{align}
&\bar{v}_{1n}= n\hat{\ell}\bar{\varphi}_n, \label{35.1}\\
&\bar{q}_n= \widehat{\mathcal{H}}_2\bar{\mathcal{H}}_{2n}+ \widehat{\mathcal{H}}_3\bar{\mathcal{H}}_{3n} -\widehat{\mathcal{E}}_1\bar{\mathcal{E}}_{1n}, \label{35.2}\\
&\bar{\mathcal{E}}_{2n}=n\big(\varepsilon s\widehat{\mathcal{H}}_3\bar{\varphi}_n - i\omega_2\widehat{\mathcal{E}}_1\bar{\varphi}_n\big) ,
\label{35.3}\\
& \bar{\mathcal{E}}_{3n}=-n\big(\varepsilon s\widehat{\mathcal{H}}_2\bar{\varphi}_n +i\omega_3\widehat{\mathcal{E}}_1\bar{\varphi}_n\big),
\label{35.4}
\end{align}
where $\hat{\ell}=s+i(\hat{v}'\cdot\omega')$. Taking into account the relation
\[
\lambda^-\bar{\mathcal{E}}_{1n}+i\omega_2\bar{\mathcal{E}}_{2n}+i\omega_3\bar{\mathcal{E}}_{3n}=0
\]
following from the last divergence constraint in \eqref{18l}, conditions \eqref{35.3} and \eqref{35.4} imply
\begin{equation}\label{35.6}
\bar{\mathcal{E}}_{1n}\lambda^-=-n\big(i\varepsilon s\hat{w}^-_\bot +\widehat{\mathcal{E}}_1|\omega'|^2\big)\bar{\varphi}_n
\end{equation}
and
\begin{equation}\label{35.7}
\widehat{\mathcal{H}}_3\bar{\mathcal{E}}_{2n}-\widehat{\mathcal{H}}_2\bar{\mathcal{E}}_{3n}=n\big(\varepsilon s |\widehat{\mathcal{H}}'|^2- i\widehat{\mathcal{E}}_1 \hat{w}^-_\bot\big)\bar{\varphi}_n,
\end{equation}
where $\hat{w}^-_\bot =\widehat{\mathcal{H}}_3\omega_2-\widehat{\mathcal{H}}_2\omega_3$.

It follows from \eqref{33}, \eqref{exp_sol_pl} and \eqref{48} that
\begin{equation}\label{50}
\left(\hat{\ell}+\frac{(\hat{w}^+)^2}{\hat{\ell}}\right)\bar{v}_{1n} -|\omega'|\bar{q}_n=0,
\end{equation}
where  $\hat{w}^+=(\widehat{H}'\cdot\omega')$. Using \eqref{35.1} and \eqref{35.2}, from \eqref{50} we derive
\begin{equation}\label{51}
n(\hat{\ell}^2+(\hat{w}^+)^2)\bar{\varphi}_{n}-|\omega'|\big(\widehat{\mathcal{H}}_2\bar{\mathcal{H}}_{2n}+ \widehat{\mathcal{H}}_3\bar{\mathcal{H}}_{3n}\big) +|\omega'|\widehat{\mathcal{E}}_1\bar{\mathcal{E}}_{1n}=0.
\end{equation}
From \eqref{34} and \eqref{exp_sol_vac} we deduce
\begin{equation}\label{57}
\varepsilon s\big(\widehat{\mathcal{H}}_2\bar{\mathcal{H}}_{2n}+ \widehat{\mathcal{H}}_3\bar{\mathcal{H}}_{3n}\big) -i\hat{w}^-_\bot \bar{\mathcal{E}}_{1n} +\lambda^-\big(\widehat{\mathcal{H}}_3\bar{\mathcal{E}}_{2n}-\widehat{\mathcal{H}}_2\bar{\mathcal{E}}_{3n}\big)=0.
\end{equation}

Combining \eqref{35.6}, \eqref{35.7}, \eqref{51} and \eqref{57} and taking into account \eqref{48}, after some algebra we obtain the equation
\[
\mathcal{L} (s,\omega')\bar{\varphi}_n=0,
\]
where
\begin{equation}\label{58}
\begin{split}
\mathcal{L} (s,\omega')= &
(\hat{\ell}^2+(\hat{w}^+)^2)\sqrt{|\omega '|^2+s^2\varepsilon^2}\\
 & +(\hat{w}^-)^2|\omega'|-\widehat{\mathcal{E}}_1^{\,2}|\omega'|^3-2i\widehat{\mathcal{E}}_1\hat{w}^-_\bot s|\omega'|\varepsilon +|\omega'|\, |\widehat{\mathcal{H}}'|^2s^2\varepsilon^2
\end{split}
\end{equation}
and $\hat{w}^-=(\widehat{\mathcal{H}}'\cdot\omega')$ (note that $(\hat{w}^-)^2=|\widehat{\mathcal{H}}'|^2|\omega'|^2-(\hat{w}^-_\bot )^2$). Using the algebraic systems for $\bar{U}_n$, $\bar{q}_n$ and $\bar{V}_n$ following from \eqref{33} and \eqref{34} as well as relations \eqref{35.1}--\eqref{35.7}, we can easily show that $\bar{\varphi}_n=0$ implies $\bar{U}_n=\bar{V}_n=0$ and $\bar{q}_n=0$. Therefore, we can construct an Hadamard-type ill-posedness example if and only if the equation $\mathcal{L} (s,\omega')=0$ has a root $s$ with $\Re\,s>0$ for some $\omega'\in\mathbb{R}^2$.

That is, the function $\mathcal{L} (s,\omega')$ in \eqref{58} is nothing else than the Lopatinski determinant, where $s$ and $\omega'$ are the Laplace and Fourier variables respectively. This function is a homogeneous function of degree three. We can thus introduce the scaling $\tilde{s}=s/|\omega'|$, $\tilde{\omega}'=\omega'/|\omega'|$. Dropping tildes, we get the equation $\mathcal{L} (s,\omega')=0$ with $|\omega'|=1$:
\begin{equation}\label{59}
\big(\hat{\ell}^2+(\hat{w}^+)^2\big)\sqrt{1+s^2\varepsilon^2}+(\hat{w}^-)^2-\widehat{\mathcal{E}}_1^{\,2}-
2i\widehat{\mathcal{E}}_1s\hat{w}^-_\bot \varepsilon + |\widehat{\mathcal{H}}'|^2s^2\varepsilon^2 =0.
\end{equation}
We seek roots of \eqref{59} in the form of a series
\[
s=s_0\varepsilon^{k_0}+s_1\varepsilon^{k_1} +s_2\varepsilon^{k_2} +\ldots
\]
for $\varepsilon \ll 1$, where $k_0<k_1<k_2<\cdots $ are rational numbers. In principle, by squaring one can reduce \eqref{59} to a polynomial equation of degree 6. This produces spurious roots of \eqref{59}, but for a polynomial equation we can apply the Newton polygon method \cite{Newton} for finding proper degrees $k_0$, $k_1$, $k_2$, ...  and then use this information for solving the original equation \eqref{59} in the form of above series. On the other hand, in our case we can manage without Newton polygons by using simple arguments. In particular, at the first stage we just rewrite equation \eqref{59} by substituting in it the Taylor series for the square root appearing there:
\begin{equation}\label{59'}
\begin{split}
  \big(s & +i(\hat{v}'\cdot\omega')\big)^2  +(\hat{w}^+)^2  +(\hat{w}^-)^2-\widehat{\mathcal{E}}_1^{\,2}  -
2i\widehat{\mathcal{E}}_1s\hat{w}^-_\bot  \varepsilon\\ &
+ s^2\left(|\widehat{\mathcal{H}}'|^2 + \Big(\big(s  +i(\hat{v}'\cdot\omega')\big)^2+(\hat{w}^+)^2\Big)\Big(\frac{1}{2}
-\frac{s^2}{8}\varepsilon^2 +\frac{s^4}{16}\varepsilon^4 -\ldots \Big)\right)\varepsilon^2  =0.
\end{split}
\end{equation}

We can easily understand that the lowest degree $k_0=0$:
\begin{equation}\label{ser}
s=s_0+s_1\varepsilon^{k_1} +s_2\varepsilon^{k_2} +\ldots \quad (0<k_1<k_2<\cdots ).
\end{equation}
Substituting \eqref{ser} into \eqref{59'} and collecting terms with zero powers of $\varepsilon$, we get the equation for $s_0$:
\begin{equation}\label{eq}
\left(s_0+i(\hat{v}'\cdot\omega')\right)^2=\widehat{\mathcal{E}}_1^{\,2}-((\hat{w}^+)^2+(\hat{w}^-)^2).
\end{equation}
This equation has a root $s_0$ with $\Re\,s_0>0$ if and only if
\begin{equation}
\widehat{\mathcal{E}}_1^{\,2}>\min_{|\omega'|=1}\left\{ (\widehat{H}'\cdot \omega')^2 + (\widehat{\mathcal{H}}'\cdot \omega')^2\right\},
\label{illcond'}
\end{equation}
where $\omega'\in\mathbb{R}^2$. Hence, condition \eqref{illcond'} is sufficient for ill-posedness. By elementary analysis (see Appendix \ref{appA}), one can show that \eqref{illcond'} is equivalent to
\begin{equation}
\widehat{\mathcal{E}}_1^{\,2}> \frac{ |\widehat{H}'|^2 +|\widehat{\mathcal{H}}'|^2-\sqrt{\big(|\widehat{H}'|^2 +|\widehat{\mathcal{H}}'|^2\big)^2-4|\widehat{H}'|^2|\widehat{\mathcal{H}}'|^2\sin^2\alpha}}{2},
\label{illcond''}
\end{equation}
where $\alpha $ is the angle formed by the vectors $\widehat{H}'$ and $\widehat{\mathcal{H}}'$. After that we easily rewrite inequality \eqref{illcond''} as \eqref{illcond} (recall that $|\widehat{H}'|^2=|\widehat{H}|^2$ and $|\widehat{\mathcal{H}}'|^2=|\widehat{\mathcal{H}}|^2$, see \eqref{constsol}). That is, condition \eqref{illcond}, which is equivalent to \eqref{illcond'}, implies ill-posedness.

Let us now
\begin{equation}
\label{e1}
\widehat{\mathcal{E}}_1^{\,2}<\min_{|\omega'|=1}\left\{ (\widehat{H}'\cdot \omega')^2 + (\widehat{\mathcal{H}}'\cdot \omega')^2\right\}.
\end{equation}
Then, for any $\omega'$ with $|\omega'|=1$ we have
\begin{equation}\label{not}
\widehat{\mathcal{E}}_1^{\,2}-((\hat{w}^+)^2+(\hat{w}^-)^2)<0.
\end{equation}
It follows from \eqref{eq} and \eqref{not} that
\[
s_0 = i\left (\pm \eta -(\hat{v}'\cdot\omega')\right):=i\tau_0^{\pm},
\]
where $\eta=\big((\hat{w}^+)^2+(\hat{w}^-)^2-\widehat{\mathcal{E}}_1^{\,2}\big)^{1/2}>0$ and $\tau_0^{\pm}\in\mathbb{R}$. Clearly, $k_1=1$ and we find $s_1$ from the equation (see \eqref{59'})
\[
\pm 2i\eta s_1+2\widehat{\mathcal{E}}_1\hat{w}^-_\bot \tau_0^{\pm}=0.
\]
Hence, $s_1=\pm i\widehat{\mathcal{E}}_1\hat{w}^-_\bot\tau_0^{\pm}/\eta:=i\tau_1^{\pm}$ and $\tau_1^{\pm}\in\mathbb{R}$.

Quite analogously we understand that $k_2=2$. Collecting terms with second powers of $\varepsilon$ in  \eqref{59'}, we obtain the equation
\[
\pm 2i\eta s_2-(\tau_1^{\pm})^2+2\widehat{\mathcal{E}}_1\hat{w}^-_\bot \tau_1^{\pm}-|\widehat{\mathcal{H}}'|^2(\tau_0^{\pm})^2-
\big((\hat{w}^+)^2-\eta^2\big)\frac{(\tau_0^{\pm})^2}{2}=0
\]
for $s_2$ which implies $s_2=i\tau_2^{\pm}$, with $\tau_2^{\pm}\in\mathbb{R}$. By finite induction we can easily show that in \eqref{ser} the degrees $k_j=j$, with $j\in \mathbb{N}$, and $s_j$ solves the equation
\[
\pm 2i\eta s_j + f_j (\tau_0^{\pm},\tau_1^{\pm},\ldots,\tau_{j-1}^{\pm})=0,
\]
where $f_j$ is a real-valued polynomial of $\tau_0^{\pm}\in\mathbb{R},\ldots , \tau_{j-1}^{\pm}\in\mathbb{R}$. Therefore, the roots of equation \eqref{59} are pure imaginary: $s=i\tau^{\pm}$ and $\tau^{\pm}\in \mathbb{R}$. That is, under condition \eqref{e1} we have no ill-posedness. Moreover, \eqref{e1} implies {\it neutral stability} of the planar plasma-vacuum interface, i.e., the Lopatinski condition holds only in a weak sense (the Lopatinski determinat has no roots $s$ with $\Re s>0$, but it has pure imaginary roots).

At last, we consider the transitional case
\begin{equation}
\label{e1'}
\widehat{\mathcal{E}}_1^{\,2}=\min_{|\omega'|=1}\left\{ (\widehat{H}'\cdot \omega')^2 + (\widehat{\mathcal{H}}'\cdot \omega')^2\right\},
\end{equation}
i.e., case \eqref{illcond=}.
If $\omega'$ is not one of the two minimum points $\omega'_*$ and $-\omega'_*$ of the function $F(\omega' )= (\widehat{H}'\cdot \omega')^2 + (\widehat{\mathcal{H}}'\cdot \omega')^2$ (see Appendix \ref{appA}), then inequality \eqref{not} holds and, as it was already proved above, we have no ill-posedness for such $\omega'$. Let $\omega' = \omega'_*$ (or $\omega' = -\omega'_*$). Then, $\widehat{\mathcal{E}}_1^{\,2}=(\hat{w}^+)^2+(\hat{w}^-)^2$ and it follows from \eqref{eq} that we have the double root $s_0=-i(\hat{v}'\cdot\omega')$.

If $(\hat{v}'\cdot\omega')=0$ (in particular, $\hat{v}'=0$), then one root of equation \eqref{59'}, which is written as
\[
  s\left\{s  - 2i\widehat{\mathcal{E}}_1\hat{w}^-_\bot  \varepsilon
+ s\left(|\widehat{\mathcal{H}}'|^2 + \big(s^2+(\hat{w}^+)^2\big)\Big(\frac{1}{2}
-\frac{s^2}{8}\varepsilon^2 +\frac{s^4}{16}\varepsilon^4 -\ldots \Big)\right)\varepsilon^2\right\}  =0,
\]
is $s=0$. By finite induction we can show that for another roots the following expansion holds:
\[
s=i\varepsilon \big( 2\widehat{\mathcal{E}}_1\hat{w}^-_\bot +\tau_2\varepsilon^2+\tau_4\varepsilon^4 +\ldots +\tau_{2k}\varepsilon^{2k}+
\ldots \big),
\]
where $\tau_{2k}\in\mathbb{R}$ for all $k\in \mathbb{N}$. Therefore, we again have neutral stability. Note that for $(\hat{v}'\cdot\omega'_*)=0$ inequality \eqref{addinst} is violated.

Let us now $(\hat{v}'\cdot\omega')\neq 0$ and we still consider case \eqref{e1'} and $\omega'=\pm \omega'_*$, i.e., $\widehat{\mathcal{E}}_1^{\,2}=(\hat{w}^+)^2+(\hat{w}^-)^2$. By analyzing \eqref{59'}, it is not difficult to understand that
\[
s=-i(\hat{v}'\cdot\omega') +s_1\sqrt{\varepsilon} +s_2\varepsilon + s_3\varepsilon\sqrt{\varepsilon} +\ldots +s_j\varepsilon^{j/2} +\ldots ,\quad j\in\mathbb{N}.
\]
For $s_1$ we obtain the equation
\[
s_1^{\,2}=2\widehat{\mathcal{E}}_1\hat{\omega}^-_{\bot}(\hat{v}'\cdot\omega').
\]
Hence, under the requirement $\widehat{\mathcal{E}}_1\hat{\omega}^-_{\bot}(\hat{v}'\cdot\omega')>0$ we have ill-posedness. In other words, for the transitional case \eqref{illcond=} the plasma-vacuum interface is violently unstable if
\begin{equation}
\widehat{\mathcal{E}}_1(\widehat{\mathcal{H}}_3\omega^*_2-\widehat{\mathcal{H}}_2\omega^*_3) (\hat{v}'\cdot\omega'_*)>0
\label{addinst}
\end{equation}
for $\omega'_*=(\omega^*_2,\omega^*_3) \in\mathbb{R}^2$ with $|\omega'_*|=1$ being the minimum point (together with $-\omega'_*$) of the function
\begin{equation}\label{F}
F(\omega' )= (\widehat{H}'\cdot \omega')^2 + (\widehat{\mathcal{H}}'\cdot \omega')^2,
\end{equation}
i.e.,
\begin{equation}\label{Fmin}
F(\omega'_*)= F(-\omega'_*)=F_{\rm min}=\min_{|\omega'|=1}\left\{ (\widehat{H}'\cdot \omega')^2 + (\widehat{\mathcal{H}}'\cdot \omega')^2\right\}.
\end{equation}
By elementary analysis (see Appendix \ref{appA}) we can exclude $\omega'_*$ from \eqref{addinst}. Using then \eqref{illcond=}, we finally get condition \eqref{illcond=1}.

Let $\widehat{\mathcal{E}}_1\hat{\omega}^-_{\bot}(\hat{v}'\cdot\omega')\leq 0$. Then
\[
s_1=\pm i \sqrt{-\widehat{\mathcal{E}}_1\hat{\omega}^-_{\bot}(\hat{v}'\cdot\omega')}.
\]
Following the above arguments, we can show that all $s_j$ for $j>1$ are also pure imaginary. This means neutral stability and completes the proof of Theorem \ref{t1}.

\begin{remark}
{\rm
If the physical restriction \eqref{physrestr} is violated, i.e., if
$\widehat{\mathcal{E}}_1^{\,2}<|\widehat{\mathcal{H}}|^2$, then the linearized problem \eqref{33}--\eqref{35} is ill-posed. Indeed,  for $\widehat{\mathcal{E}}_1^{\,2}>|\widehat{\mathcal{H}}|^2$, one has:
\[
\widehat{\mathcal{E}}_1^{\,2}>|\widehat{\mathcal{H}}'|^2\geq (\widehat{\mathcal{H}}'\cdot \omega'_*)^2=F(\omega'_*)\geq F_{\rm min},
\]
where $\omega'_*\in \mathbb{R}^2$ is such that $\widehat{H}'\cdot \omega'_*=0$ and $|\omega'_*|=1$, and  $F(\omega' )$ and $F_{\rm min}$ were defined in \eqref{F} and \eqref{Fmin}. In fact,  we have ill-posedness even if $\widehat{\mathcal{E}}_1^{\,2}=|\widehat{\mathcal{H}}|^2$ because for this case it follows from the condition $2\hat{p} +|\widehat{H}|^2=|\widehat{\mathcal{H}}|^2-\widehat{\mathcal{E}}_1^{\,2}$ that $\widehat{H}=0$ (and $\hat{p}=0$, where $\hat{p}$ is a constant pressure). For $\widehat{H}=0$ the ill-posedness condition \eqref{illcond} is reduced to the true inequality $\widehat{\mathcal{E}}_1^{\,2}>0$ for any nonzero $\widehat{\mathcal{E}}_1$.}
\label{r4}
\end{remark}

\section{Concluding remarks and open problems}

\label{s5}

The ill-posedness condition \eqref{illcond} is written in such a form that gives us the following guess about the well-posedness condition for the nonlinear problem \eqref{22}--\eqref{24}:
\begin{equation}\label{well_nonlin}
{\mathcal{E}}_1^{\,2}< \frac{ |{H}|^2 +|{\mathcal{H}}|^2-\sqrt{\big(|{H}|^2 +|{\mathcal{H}}|^2\big)^2-4|{H}\times{\mathcal{H}}|^2}}{2}\quad \mbox{at}\ x_1=0.
\end{equation}
Our hypothesis is that the nonlinear problem is well-posed (locally in time) provided that the initial data satisfy the ``stability'' condition \eqref{well_nonlin} (together with other necessary conditions like compatibility conditions, regularity assumptions, etc.). Note  that if ${H}\times{\mathcal{H}}=0$ at some point of a non-planar interface, then our hypothetical necessary and sufficient well-posedness condition \eqref{well_nonlin} is violated.

The equality \eqref{illcond=} just defines a codimension-one set in the space of seven parameters $(\widehat{H}',\widehat{\mathcal{H}}',\hat{v}',\widehat{\mathcal{E}}_1)\in\mathbb{R}^7$ of the linear problem \eqref{33}--\eqref{35}. In this sense the study of the transitional case \eqref{illcond=} is not important for possible future nonlinear analysis under the (hypothetical) well-posedness condition \eqref{well_nonlin} (clearly, the counterpart of \eqref{well_nonlin} with non-strict inequality is a bad assumption for initial data of the nonlinear problem to be satisfied for $t>0$).

There are different possible ways to prove the local-in-time well-posedness of the nonlinear plasma-vacuum problem under the ``stability'' condition \eqref{well_nonlin}. In our opinion, one can try to use approaches applied in \cite{CDS2,CS1,CS2,GuWang,MTT1,SWZ1,SWZ2}. We think that this is a very interesting challenging open problem for future research. Even the proof of an a priori estimate for the constant coefficient problem \eqref{33}--\eqref{35} under condition \eqref{well_nonlin} satisfied for the constant solution is a difficult problem, in particular, because of the fact that system \eqref{33} is not hyperbolic and, hence, Kreiss' symmetrizers technique \cite{Kreiss} is not directly applicable. At the same time, it is clear how to prove well-posedness for initial data with a sufficiently small $|\mathcal{E}_{1}|$ at $x_1=0$ (and under the non-collinearity condition \eqref{noncollin} at $t=0$). For this purpose, we can follow the strategy of \cite{CDS,M1,MT} combined with the idea of a slightly compressible (isentropic) regularization \cite{MTTpv} of the linearized problem. In our opinion, this is not so interesting in comparison with the main open problem described above.

\section*{Acknowledgments} This work was supported by RFBR (Russian Foundation for Basic Research) grant No. 19-01-00261-a.

\appendix
\section{Exclusion of the wave vector from the violent instability conditions}

\label{appA}

Let $F(\omega' )= (\widehat{H}'\cdot \omega')^2 + (\widehat{\mathcal{H}}'\cdot \omega')^2$. Our goal is to find
\[
F_{\min}=\min\limits_{|\omega'|=1}F(\omega' ),
\]
cf. \eqref{illcond'}. Let $\alpha \in [0,2\pi ]$ be the angle formed by the vectors $\widehat{H}'$ and $\widehat{\mathcal{H}}'$ and $x$ be the angle formed by the bisector of $\alpha$ and the wave vector $\omega'$. Then for $|\omega'|=1$ we have
\[
F(\omega' )=f(x)=a\cos^2\left(x+\frac{\alpha}{2}\right)+b\cos^2\left(x-\frac{\alpha}{2}\right) = \frac{a+b}{2}+\frac{
a\cos (2x+\alpha)+b\cos (2x -\alpha )}{2},
\]
where $a=|\widehat{H}'|^2$ and $b =|\widehat{\mathcal{H}}'|^2$. Clearly, the trigonometrical function $f(x)$ has two minimum points on the interval for $[0,2\pi ]$ (one of them corresponds to some $\omega'$ and the another  one to $-\omega'$).

We seek extreme points of $f(x)$:
\[
-f'(x)=a\sin (2x+\alpha)+b\sin (2x -\alpha )=0.
\]
We easily rewrite this equation as
\[
\tan 2x= \frac{b-a}{b+a}\tan \alpha.
\]
Let $x_*\in[0,2\pi ]$ be a solution of the last equation. One can show that for both minimum points
\begin{equation}
\cos 2x_*=-\frac{(b+a)\cos\alpha}{\sqrt{(a+b)^2 -4ab\sin^2\alpha}}\quad \mbox{and}\quad
\sin 2x_*=-\frac{(b-a)\sin\alpha}{\sqrt{(a+b)^2 -4ab\sin^2\alpha}}.
\label{A1}
\end{equation}
Using \eqref{A1}, we find
\[
\begin{split}
F_{\min} & =\frac{a+b}{2} +\frac{(a+b) \cos 2x_*\cos\alpha+(b-a)\sin 2x_*\sin\alpha}{2} \\[6pt] &=\frac{a +b -\sqrt{(a+b)^2 -4ab\sin^2\alpha}}{2}.
\end{split}
\]
We have thus obtained the violent instability condition \eqref{illcond''} which is rewritten as \eqref{illcond}.

Let condition \eqref{illcond=} be satisfied, i.e.,
\begin{equation}
\widehat{\mathcal{E}}_1^{\,2}=\frac{a +b -\sqrt{(a+b)^2 -4ab\sin^2\alpha}}{2}.
\label{A2}
\end{equation}
We rewrite inequality \eqref{addinst} for the minimum points  as
\[
\widehat{\mathcal{E}}_1|\widehat{\mathcal{H}}'|\,|\hat{v}'|\sin\left(\frac{\alpha}{2}-x_*\right)\cos\left(x_*-\frac{\alpha}{2}-\beta\right)=
\frac{1}{2}\widehat{\mathcal{E}}_1|\widehat{\mathcal{H}}'|\,|\hat{v}'|(\sin(\alpha +\beta-2x_*)-\sin\beta )>0,
\]
where $\beta$ is the angle formed by the vectors $\hat{v}'$ and $\omega'$, and $x_*$ is such that \eqref{A1} holds. By virtue \eqref{A1} and \eqref{A2}, the last inequality can be further reduced to
\[
\widehat{\mathcal{E}}_1|\widehat{\mathcal{H}}'|\,|\hat{v}'|\big((\widehat{\mathcal{E}}_1^{\,2}-b)\sin\beta -a\sin (\alpha +\beta )\cos\alpha  \big) >0
\]
that is rewritten as
\begin{equation}\label{illcond=1'}
\widehat{\mathcal{E}}_1 \big( (\widehat{H}'\cdot\widehat{\mathcal{H}}')(\hat{v}_2\widehat{H}_3-\hat{v}_3\widehat{H}_2)+
(|\widehat{\mathcal{H}}|^2-\widehat{\mathcal{E}}_1^{\,2})
(\hat{v}_2\widehat{\mathcal{H}}_3-\hat{v}_3\widehat{\mathcal{H}}_2)\big) >0,
\end{equation}
Inequality \eqref{illcond=1'} is nothing else than  condition \eqref{illcond=1}.

\end{document}